# Intransitive Machines

Alexander Poddiakov

Department of Psychology
National Research University Higher School of Economics

**Abstract**
The intransitive cycle of superiority is characterized by such binary relations between A, B, and C that A is superior to B, B is superior to C, and C is superior to A (i.e., A>B>C>A—in contrast with transitive relations A>B>C). The first part of the article presents a brief review of studies of intransitive cycles in various disciplines (mathematics, biology, sociology, logical games, decision theory, etc.), and their reflections in educational materials. The second part of the article introduces the issue of intransitivity in elementary physics. We present principles of building mechanical intransitive devices in correspondence with the structure of the Condorcet paradox, and describe five intransitive devices: intransitive gears; levers; pulleys, wheels, and axles; wedges; inclined planes. Each of the mechanisms are constructed as compositions of simple machines and show paradoxical intransitivity of relations such as "to rotate faster than", "to lift", "to be stronger than" in some geometrical constructions. The article is an invitation to develop teaching materials and problems advancing the understanding of transitivity and intransitivity in various areas, including physics education.

**Key words:** intransitivity of superiority; rock-paper-scissors game; simple machines; intransitive machines; physics education

**Introduction: Intransitivity of Superiority as an Interdisciplinary Issue**

The intransitive cycle (or loop) of superiority is characterized by such binary relations between A, B, and C that A is superior to B, B is superior to C, and C is superior to A (i.e., A>B>C>A—in contrast with transitive relations A>B>C). Although some examples of the intransitivity of superiority have long been described, the phenomenon really came into light and became the focus of interdisciplinary discussion as a result of articles published by Martin Gardner, one of the most famous popular mathematics writers. In 1970 and 1974, he published two articles about paradoxical objects violating transitivity in his "Mathematical Games" column in *Scientific American* (Gardner, 1970, 1974). The articles were reproduced in "The Colossal Book of Mathematics" (Gardner, 2001). The mathematical facts described in the articles changed many people's perceptions about the world: to their amazement, they learned that in sets of some objects (Efron's dice, playing cards, roulette wheels, voters' preferences, etc.) elements of a set beat one another in a counterintuitive way — like in the rock-paper-scissors game. For example, in a set of dice A, B, and C, Die A more often rolls a higher number than Die B in the pair A-B, Die B more often rolls a higher number than Die C in the pair B-C, but Die C more often rolls a higher number than Die A in the pair A-C.

Perhaps the most famous — and relatively simple — version of combinations of numbers on intransitive dice was designed by Bradley Efron. Efron's dice are presented in the National Museum of Mathematics (USA). On the Microsoft site *Non-Transitive Dice* (2011), these are cited to highlight how important it is to analyze intransitivity by means of artificial intelligence in large sets of data.

Blue dice: 4, 4, 4, 4, 0, 0
Yellow dice: 3, 3, 3, 3, 3, 3
Red dice: 6, 6, 2, 2, 2, 2
Green dice: 5, 5, 5, 1, 1, 1

Above one can see that the probability of the blue die winning when paired with the yellow die is 2/3 (the blue die rolls a higher number more often than the yellow die). The same goes for the yellow die paired with the red die and the red die paired with the green die. When the green die is pitted against the blue die, the odds are also 2/3. Thus, the blue die beats the yellow die, the yellow die beats the red die, the red die beats the green die, and the green die beats the blue die. This is an example of the intransitive cycle of superiority.

The number of intransitive mathematical objects is infinite. For example, one can design cyclic chains of intransitive multi-sided dice of arbitrary length: "It is possible to find $n$ dice with $n$ faces each such that these dice form a set of intransitive dice" (Deshpande, 2000, p. 5). Cyclic chains of five and more intransitive dice can include subsets of shorter cyclic chains (e.g., of three elements), and it is possible that within the subset the direction in which elements "beat" each other may be opposite to that within the whole chain. It may also happen within such sets of intransitive dice that their doubling (situations when Player 1 rolls not the only die A but two dice A at the same time, Player 2 rolls two dice B at the same time, etc.) leads to a reversal in the direction of beating: A beats B, B beats C, C beats D, and D beats A, but the double set AA loses to BB, BB loses to CC, CC loses to DD, and DD loses to AA (Grime, 2017).

Statistical intransitive relations can be transformed into deterministic intransitive relations — for example, between lengths of sticks. Many math educators (and not only they) know standard tests on transitive inference that require the learner to make the correct conclusion from premises "Stick A is longer than B, and Stick B is longer than C" in order to answer the question "Which of the sticks is the longest one?" Certainly, transitivity does work here in a comparison of three sticks. Yet intransitivity emerges if we consider more complex situations such as those comparing sets that each include three and more elements. Three or more sets of three or more sticks (e.g., pencils) whose lengths are equal to the numbers on the faces of intransitive dice (or numbers on intransitive playing cards, or numbers in cells of some magic squares, etc.) will be in deterministic intransitive relations of superiority of length. If one organizes a tournament between "teams" of pencils and compares the length of each pencil from each set with the length of each pencil from the other sets, s/he can see that pencils of Set A are more often longer than pencils of Set B, which are more often longer than pencils of Set C, and the latter tend to be longer than those of Set A (see Figure 1). In other words, the results of the tournament between such pencil teams will be intransitive: Team A beats B, B beats C, and C beats A. Thus, the elementary relation "to be longer than" is transitive, but the more complex relation "to be longer more often than" can be intransitive.

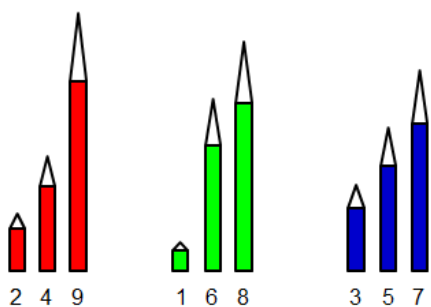

**Fig. 1** An example of intransitive sets of sticks (e.g., pencils)
Numbers to define the pencils' lengths are taken from the magic square presented by Gardner (1974, p. 120)

A note on terminology: the terms "intransitive" and "non-transitive" (e.g., "intransitive dice" and "non-transitive dice") are used for such sets as synonyms in the math literature in spite of some difference between the logical terms "intransitive relation" and "non-transitive relation". "An intransitive relation is one that does not hold between *a* and *c* if it also holds between *a* and *b* and between *b* and *c* for any substitution of objects for *a*, *b*, and *c*. ... A non-transitive relation is one that may or may not hold between *a* and *c* if it also holds between *a* and *b* and between *b* and *c*, depending on the objects substituted for *a*, *b*, and *c*. In other words, there is at least one substitution in which the relation between *a* and *c* does hold and at least one substitution in which it does not" (Transitive law). Examples of intransitive relations are "to be the (biological) daughter of" and "is the square of", and examples of non-transitive relations are "… loves …" and "is not equal to" (Ibid). In this article we will use the term "intransitive" as explicitly related to the concept of intransitive cycles.

Mathematical studies of intransitive objects are actively developing. They show that the issue of intransitivity holds opportunities for new discoveries in studies of large sets of dice, including multi-sided dice (Buhler et al., 2018; Conrey et al., 2016), and opportunities of posing interesting and novel mathematical problems which may reveal unknown and important mathematical relationships (see the discussion of this issue by Fields medalist T. Gower (2017) on the pages of his Polymath project).

Intransitivity is actively studied not only in mathematics, but also in other disciplines. In theoretical physics, an area closely related to mathematics, advanced models of intransitive relations are developed to explain complex processes of thermodynamics (Klimenko, 2013). In biology, the phenomenon of intransitive competition between biological species and between individuals is studied at different levels of biological and ecological interactions, and the metaphorical principle of the rock-paper-scissors game is viewed within this paradigm as a crucial condition of biodiversity (Allesina & Levine, 2011; Kerr et al., 2002; Kirkup & Riley, 2004; Jiang et al., 2011; Maynard et al., 2017; Reichenbach et al., 2007; Sinervo & Lively, 1996; Szolnoki et al., 2014; Vandermeer, 2011). Interestingly, "robot Darwinism", a phenomenon analogical to biological intransitive competition, is observed (but not yet scientifically studied) in fights of remote-controlled vehicles in BattleBots shows. More often than not, "lifter" vehicles tend to beat "spinners", "spinners" tend to beat "pure wedges", and "pure wedges" often beat "lifters" (Atherton, 2013).

The intransitivity between the participants of sporting events is well known. For instance, at a chess tournament, Player A might tend to win when playing with Player B, B might tend to beat C, while C tends to beat A (West & Hankin, 2008). However, in logical board games such as chess or checkers, the concept of intransitivity also applies to positions. It works in the following way: Position A for White is preferable to Position B for Black (i.e., when offered a choice, one should choose A), Position B for Black is preferable to Position C for White, which is preferable to Position D (Black) – but the latter is preferable to Position A (White) (Poddiakov, 2017). In chess, the number of such intransitive chains is extremely high (on a par with the number of all possible positions on the chess board) and their length can vary greatly, from a four-position chain involving only a king and a pawn for White and Black to astronomically long chains with a number of pieces on both sides (Filatov, 2017).

The significance of objective intransitivity of superiority for decision theory is discussed by Anand (1993), Bar-Hillell and Margalit (1988), Fishburn (1991), and Temkin (1996, 2012). Fishburn sees an analogy between our current understanding of the limits imposed by the formal-logical axiom of transitivity applied to fairly simple objects and phenomena and the move from Euclidean to non-Euclidean geometry. According to him, this is similar to the move away from viewing the Newtonian worldview as absolute (Fishburn 1991). Interdisciplinary reviews and analyses of intransitive relations in various areas (from math to biology, sociology, and psychology) are presented by Fisher (2008), Klimenko (2014, 2015), Poddiakov & Valsiner (2013), and West & Hankin (2008).

Concerning education, perhaps T. Roberts presents the most compelling reason for teaching intransitivity: "Transitivity and intransitivity are fascinating concepts that relate both to mathematics and to the real world we live in", and teachers can "interest and engage students of almost any age, as they seek to discover which relationships are transitive, and which are not, and further to try to discover any general rules that might distinguish between the two" (Roberts 2004, p. 63).

Various problems ranging in complexity are designed to promote intransitivity understanding in students. Relatively simple math problems on intransitivity can be presented to secondary school students (Beardon, 1999/2011; Scheinerman, 2012). An advanced example of a complex problem based on the knowledge from different domains (biology, physics, and high school mathematics) is presented by Strogatz in his book *Nonlinear Dynamics and Chaos: With Applications to Physics, Biology, Chemistry, and Engineering*. The beginning of the problem is formulated as an interdisciplinary introduction:

"In the children's hand game of rock-paper-scissors, rock beats scissors (by smashing it); scissors beat paper (by cutting it); and paper beats rock (by covering it). In a biological setting, analogs of this non-transitive competition occur among certain types of bacteria (Kirkup & Riley, 2004) and lizards (Sinervo & Lively, 1996). Consider the following idealized model for three competing species locked in a life-and-death game of rock-paper-scissors…" (Strogatz, 2015, pp. 191-192).

Besides advanced biological models, thermodynamics, statistics and high school mathematics, what about elementary physics? Below we describe mechanisms in intransitive relations like objects in the rock-paper-scissors game, intransitive dice and other intransitive games analyzed by logicians and mathematicians. In contrast with intransitive mathematical objects which show the intransitivity of relations such as "to roll a higher number more often", intransitive machines show intransitivity of relations such as "to rotate faster than", "to lift", and "to be stronger than".

All the machines described below are constructed as compositions of simple machines: levers, wheels and axles, pulleys, wedges, and inclined planes. Mastery of simple machines is considered important for people of various ages and professions: from preschool, primary, and secondary school children to university students and professionals whether in traditional industries or newer fields like biomechanics, materials science, and design of AI systems recognizing images (Doudna, 2015; Field & Solie, 2007; Korur et al., 2017; Marulcu & Barnett, 2013; Okuno & Fratin, 2014; Olteţeanu, 2016; Rogers, 2005; Wakild, 2006; Williams et al., 2003).

Intransitive mechanisms show the possibilities of simple machines by offering an unexpected view of possible geometrical combinations of the machines and their elements. They are tools which can be used to deepen understanding of geometry and mechanics and to show students that compositions of familiar devices can work in counter-intuitive ways. They can also help to advance students' scientific beliefs about intransitive relations.

**Condorcet-like Compositions**

All the mechanisms described below are constructed as Condorcet-like compositions, in correspondence with the structure of elements in the Condorcet paradox (or the voting paradox). A description of the Condorcet paradox which can be understood by secondary school students is presented by Beardon (1999/2011):

"In the … Voting Paradox there are 3 candidates for election. The voters have to rank them in order of preference. Consider the case where 3 voters cast the following votes: ABC, BCA and CAB:
A beats B by 2 choices to 1
B beats C by 2 choices to 1
but A cannot be the preferred candidate because A loses to C, again by 2 choices to 1.
This is an example of intransitivity".

A body of work has been devoted to the voting paradox and its meaning for logic, social sciences, and practices (see a review and analysis in (Gehrlein, 2006)). Here we will show the application of its structure for the design of intransitive mechanisms. We will use chains of mechanical elements ordered in correspondence with the Condorcet paradox:

ABC
BCA
CAB

One can see that the first element of any line goes to the last position in the next line and moves all the other elements one position to the right without changing the sequence. Let us apply it to mechanics.

**Intransitive Gears (or Friction Wheels)**

Let us consider three double gears (or friction wheels), using the following notation of elements:
X is a larger gear (a larger wheel),
Y is a smaller gear (a smaller wheel),
Z is an empty part of a shaft (without any gear or wheel on it).

Then, in correspondence with the Condorcet paradox:
the first double-gear (A) will have the element sequence X, Y, Z;
the second double-gear (B) will have the element sequence Z, X, Y; and
the third double-gear (C) will have the element sequence Y, Z, X.

Figure 2 shows that, if joined in pairs, A's rotational speed is higher than B's in the pair A-B; the rotational speed of B is higher than that of C in the pair B-C; but the rotational speed of C is higher than rotational speed of A in pair A-C (Poddiakov, 2010; Poddiakov &Valsiner, 2013). Here, as in other intransitive games (rock-paper-scissors, intransitive dice, lotteries, etc.), it is preferable to choose A in the pair A-B, B in the pair B-C, and C in the pair A-C (if the "winner" is the fastest object in a pair). Yet in contrast with games such as intransitive dice and lotteries, our system is mechanical and deterministic. To win in this game, one has to comprehend geometry and mechanics.

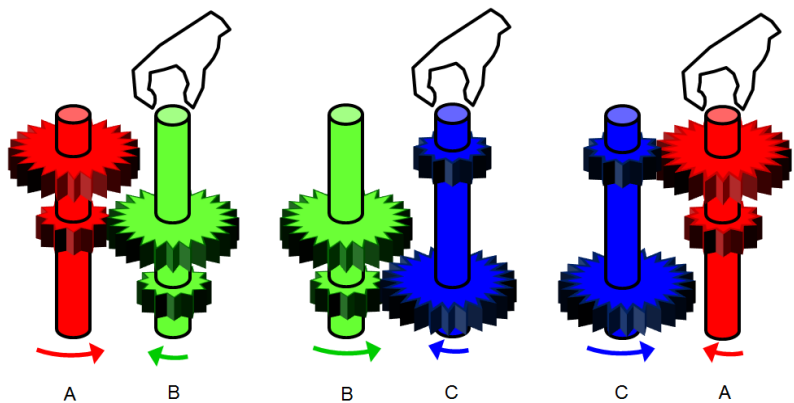

**Fig. 2** Intransitive gears
Gear A rotates faster than B in pair A-B, B rotates faster than C in pair B-C, and C rotates faster than A in pair A-C

It is important to note that joining all three double-gears (or friction wheels) will lead to a jam. Intransitivity is a property of binary relations, so we test it in pairs. Similarly, intransitivity is not revealed by rolling all intransitive dice belonging to one set at once, or by all players throwing all intransitive playing cards on the table at once.

**Intransitive Double Levers**

The way gear trains work is based on the law of the lever; a gear can be considered as a lever. Let us transform intransitive double gears into intransitive double levers using the same Condorcet-like structure.
Each of the double levers has a shaft whose long end is to be inserted into a hole and two levers of different lengths (see Figure 3). When the levers are joined two by two, given that they are equal in torque created manually by participants, the double Lever A will be stronger than B (and will move it in spite of counteraction), the double Lever B will be stronger than C, and C will be stronger than A. Again, if playing, it is preferable to let the other player choose the lever. The relation between the lengths of the double levers' arms can be such that whichever lever the strongest member of one of the teams chooses, the weakest member of the other team can nevertheless select a more advantageous double gear and win.

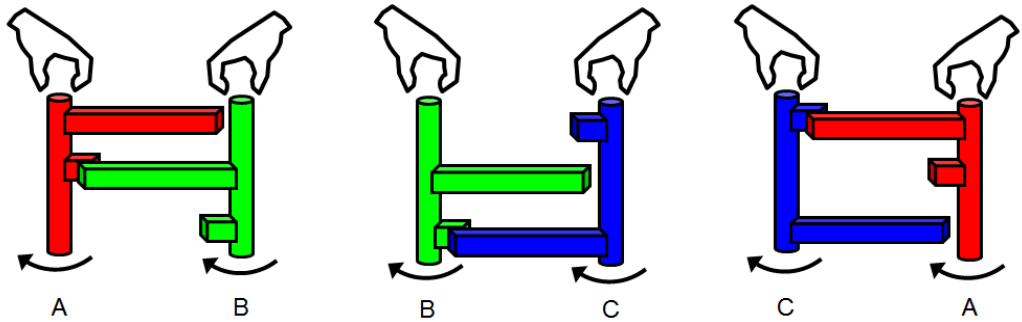

**Fig. 3** Intransitive double levers
With the same rotation force applied to the shaft, Lever A will overpower Lever B as A's lever arm (the perpendicular distance from the fulcrum to the line of action of the effort) is shorter than B's lever arm. Lever B will overpower Lever C and Lever C will overpower Lever A

**Intransitive Machines with Pulleys, Wheels, and Axles**

In the example above, the levers were turned by hand. This work can also be carried out by the force of gravity with the help of pulleys.
   To this effect, we can add the shafts of the intransitive double gears with the identical pulleys and loads on cables. We then insert the shafts of the intransitive wheels and axles into holes in a vertical wall (see Figure 4). In a way it looks like elements of a mechanical wall clock, but this is not the case. One can see that Load A will lift Load B in the pair A-B ("A is stronger than B"), Load B will lift Load C in the pair B-C ("B is stronger than C"), and Load C will lift Load A in the pair A-C ("C is stronger than A").

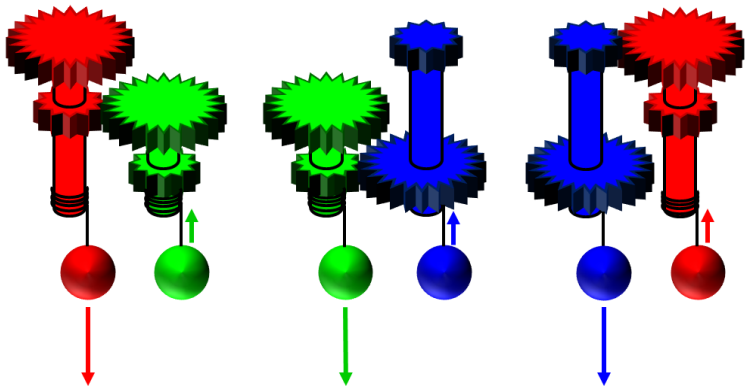

**Fig. 4** Intransitive machines with pulleys, wheels and axles
In each pair of double wheels, the weight on the shaft where the gear train involves the smaller gear (the one with fewer teeth) will go down faster and further compared to the weight going up on the shaft where the gear train involves the bigger gear (the one with more teeth)

   A frequently asked question concerning the composition of these pulleys, wheels, and axles has to do with the possibility of constructing a perpetual motion machine ("If B lifts A, C lifts B, and A lifts C, can we use it to make a perpetual motion machine?"). The answer "No, we cannot" is based on the law of the lever. In each pair (A-B, B-C, and A-C), the mechanical advantage of one member of the pair is related to the greater distance of lowering its own load and the lesser distance of lifting the other load of the same weight. Thus, the system's center of gravity moves down, and its total potential energy decreases in each pair. In the end, all the loads will be at the lower point (e.g., on the floor), similar to the weights of mechanical wall clocks (in this regard, the clocks and the pulleys-wheels-axles are similar). So, the example of the intransitive "pulleys-wheels-axles" is of such local (in pairs) mechanical advantages and proportional losses in distances in which the law of the lever works in a seemingly counterintuitive way because of the geometrical design of the devices. If one wishes, it can be used to explain laws of conservation along with other examples.

**Intransitive Machines with Wedges**

Let us consider three stylized Mobile Assault Towers having long wedges as weapons (something like sharpened rams or spears), such as the ones in Figure 5; these models were inspired by pictures of Assyrian battering rams. One can see that the spear of Tower A can touch Tower B but not vice versa, the spear of Tower B can touch Tower C but not vice versa, and the spear of Tower C can touch Tower A but not vice versa. Thus, A beats B, B beats C, and C beats A (Poddiakov, 2010; Poddiakov & Valsiner, 2013).

For demonstration, one can make Mobile Assault Towers from plastic or wooden geometrical profiles with inserted felt-tip pens of different colors. Tower A will mark Tower B but not vice versa, Tower B will mark Tower C but not vice versa, and Tower C will mark Tower A but not vice versa.

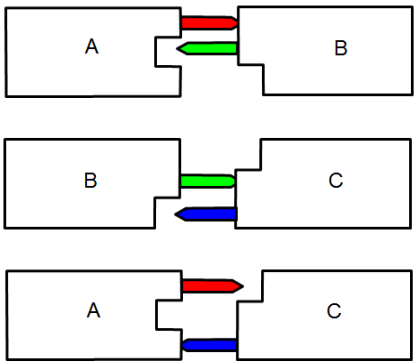

**Fig. 5** Intransitive machines with wedges (sharpened rams or spears)

**Intransitive Ramps (Inclined Planes)**

Let us consider three stylized combs with some teeth broken off and the remaining teeth having either a rectangular or wedge-like shape (see Figure 6). One can see that, in the case of a frontal collision, Comb A can serve as a ramp for Comb B and lift it (but not vice versa), Comb B can serve as a ramp for Comb C and lift it (but not vice versa), and Comb C can serve as a ramp for Comb A and lift it (but not vice versa).

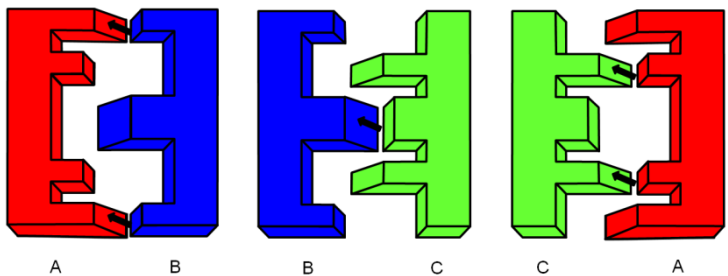

**Fig. 6** Intransitive ramps (inclined planes)

**Conclusion**

All the intransitive machines presented above are constructed as Condorcet-like compositions. They are intransitive in the same sense that objects in various mathematical intransitive games (intransitive dice, intransitive playing cards, etc.) are intransitive, but are constructed as deterministic mechanical devices. Based on a chosen criterion (e.g., "to rotate faster", "to lift and be not lifted", "to achieve more mechanical advantage", etc.), one should prefer A in pair A-B, B in pair B-C, and C in pair A-C. Of course, this set of intransitive mechanisms can be supplemented with others that are perhaps more paradoxical and counterintuitive.

Students can carry out various activities on the topic of intransitive mechanical devices and mechanical intransitivity. As they study intransitive mechanical devices or look at pictures they can consider how the devices might work in an intransitive way, establish which sequence of pairs form an intransitive loop, and figure out the direction of "advantage" in the loops. Students can also invent their own devices, modified or original as the case may be.

For example, based on the knowledge of intransitive cycles of domination in a group of stylized birds (Individual A dominates Individual B, B dominates C, and C dominates A), one can construct three mechanical models. The birds stand on the ground with open wings. In the case of a frontal meeting, two birds touch one another with their wings and then the following "intransitivity of domination" is observed: Bird B pulls its head in front of Bird A but not vice versa, Bird C pulls its head in front of Bird B but not vice versa, and Bird A pulls its head in front of Bird C but not vice versa. A tip: one wing of each bird is immovable relative to its body, while the other wing can move and is mechanically joined with the bird's head. The locations of where the points of the birds' wings touch one another are defined in correspondence with the Condorcet structure.

Some problems for students can be based on problems formulated for intransitive dice. For example, one can try to invent a chain of more than three intransitive double gears (levers), perhaps four or five, so that the last element rotates faster (is stronger) than the first one, and to invent an algorithm of construction of such a chain of arbitrary length.

The problem of designing multiple player games has already been posed for intransitive dice (Bozóki, 2014, p. 39). M. Oskar van Deventer invented "a three player game, a set of dice where two of your friends may pick a die each, yet you can always pick a die that has a better chance of beating both opponents at the same time" (Grime, 2017, p. 4). Such sets have also been designed that "for any three players, there is a fourth that beats all of them" (Pegg, 2005). Respectively, one can attempt to invent sets of intransitive mechanical devices (levers, wheel and axles, stylized birds, etc.) such that two (or three) players may pick a device each, yet the third (or fourth) one can always pick a device that beats (dominates) all of their opponents.

There are also general problems to be discussed with students. Are Condorcet-like compositions the only ones which provide an opportunity to construct intransitive mechanical and geometrical constructions? Are other principles possible?

What about mechanical intransitivity at other structural levels of interactions? For example, what about intransitivity of hardness of composite materials such that Material A abrades B rather than B abrades A, B abrades C rather than C abrades B, and C abrades A rather than A abrades C? We mean here not intransitivity of resistance caused by different qualities of making such materials, for example, at three different factories (Trybuła, 1961) but intransitivity of hardness as a consequence of the inherent properties of some composite materials. Is it possible or not? Why? It seems that answers to these questions would require knowledge of not only mechanics, but also chemistry, high school mathematics, and other components of materials science. This is by no means an exhaustive list of problems.

Our article is an invitation to more deeply explore the issue of intransitivity in education, including physics education, and also to develop teaching materials and problems advancing the understanding of transitivity and intransitivity in various disciplines.

**References**

1. Allesina, S., & Levine, J. M. (2011). A competitive network theory of species diversity. *Proceedings of the National Academy of Sciences*, *108*(14), 5638-5642. https://doi.org/10.1073/pnas.1014428108
2. Atherton, K. D. (2013). A brief history of the demise of battle bots. Retrieved from http://www.popsci.com/technology/article/2013-06/elaborate-history-how-wedges-ruined-battlebots
3. Bar-Hillel, M., & Margalit, A. (1988). How vicious are cycles of intransitive choice? *Theory and decision*, *24*(2), 119-145. https://doi.org/10.1007/BF00132458.
4. Beardon, T. (1999/2011). Transitivity. Retrieved from http://nrich.maths.org/1345.
5. Bozóki, S. (2014). Nontransitive dice sets realizing the Paley tournaments for solving Schütte's tournament problem. *Miskolc Mathematical Notes, 15*(1), 39–50. Retrieved from http://real.mtak.hu/62873.
6. Buhler, J., Graham, R., & Hales A. (2018). Maximally nontransitive dice. *The American Mathematical Monthly*, *125*(5), 387-399. https://doi.org/10.1080/00029890.2018.1427392.
7. Conrey, B., Gabbard, J., Grant, K., Liu, A., & Morrison K. (2016). Intransitive Dice. *Mathematics Magazine*, *89*(2), 133-143. https://doi.org/10.4169/math.mag.89.2.133.
8. Deshpande, M. N. (2000). Intransitive dice. *Teaching statistics*, 22 (1), 4–5. https://doi.org/10.1111/1467-9639.00002.
9. Doudna, K. (2015). *The Kids' Book of Simple Machines: Cool Projects & Activities that Make Science Fun!* Minneapolis, MN: Mighty Media Kids.
10. Field, H. L. & Solie, J. B. (2007). *Introduction to Agricultural Engineering Technology: A Problem Solving Approach.* New York: Springer Science+Business Media.
11. Filatov, A. (2017). [Intransitive chess positions.] *Nauka i Zhizn', 7*, 117-120. Retrieved from https://www.nkj.ru/archive/articles/31727. (in Russian).
12. Fishburn, P. C. (1991). Nontransitive preferences in decision theory. *Journal of risk and uncertainty*, *4*(2), 113–134. https://doi.org/10.1007/BF00056121.
13. Fisher, L. (2008). *Rock, Paper, Scissors: Game Theory in Everyday Life*. New York: Basic books.
14. Gardner, M. (1970). The paradox of the nontransitive dice and the elusive principle of indifference. *Scientific American*, *223*(6), 110-114.
15. Gardner, M. (1974). On the paradoxical situations that arise from nontransitive relations. *Scientific American*, *231*(4), 120-125.
16. Gardner, M. (2001). *The Colossal Book of Mathematics*. New York: W.W. Norton.
17. Gehrlein, W. V. (2006). *Condorcet's Paradox*. Berlin: Springer-Verlag Berlin Heidelberg.
18. Gower, T. (2017). A potential new Polymath project: intransitive dice. Retrieved from https://gowers.wordpress.com/2017/04/28/a-potential-new-polymath-project-intransitive-dice.
19. Grime, J. (2017). The bizarre world of nontransitive dice: games for two or more players. *The College Mathematics Journal, 48*(1), 2-9. https://doi.org/10.4169/college.math.j.48.1.2.
20. Jiang, L.-L., Zhou, T., Perc, M., & Wang, B.-H. (2011). Effects of competition on pattern formation in the rock-paper-scissors game. *Physical Review E*, *84*(2). https://link.aps.org/doi/10.1103/PhysRevE.84.021912.
21. Kerr B., Riley, M. A., Feldman M. W., & Bohannan B. J. M. (2002). Local dispersal promotes biodiversity in a real-life game of rock–paper–scissors. *Nature, 418*, 171-174. http://dx.doi.org/10.1038/nature00823.
22. Kirkup, B. C., & Riley, M. A. (2004). Antibiotic-mediated antagonism leads to a bacterial game of rock–paper–scissors *in vivo*. *Nature, 428*, 412-414. http://dx.doi.org/10.1038/nature02429.
23. Klimenko, A. Y. (2013). Complex competitive systems and competitive thermodynamics. *Philosophical Transactions of Royal Society A, 371*: 20120244. https://doi.org/10.1098/rsta.2012.0244.

24. Klimenko, A. Y. (2014). Complexity and intransitivity in technological development. *Journal of Systems Science and Systems Engineering, 23*(2), 128-152. https://doi.org/10.1007/s11518-014-5245-x.
25. Klimenko, A. Y. (2015). Intransitivity in theory and in the real world. *Entropy, 17*(6), 4364-4412. https://doi.org/10.3390/e17064364.
26. Korur, F., Efe, G., Erdogan, F., & Tunç, B. (2017). Effects of toy crane design-based learning on simple machines. *International Journal of Science and Mathematics Education, 15*(2), 251–271. https://doi.org/10.1007/s10763-015-9688-4.
27. Marulcu, I., & Barnett, M. (2013). Fifth graders' learning about simple machines through engineering design-based instruction using LEGO™ materials. *Research in Science Education*, *43*(5), 1825–1850. https://doi.org/10.1007/s11165-012-9335-9.
28. Maynard, D. S., Crowther, T. W., & Bradford, M. A. (2017). Competitive network determines the direction of the diversity-function relationship. *Proceedings of the National Academy of Sciences of the United States of America, 114*(43), 11464–11469. http://dx.doi.org/10.1073/pnas.1712211114.
29. Non-transitive dice. (2011). Retrieved from https://www.microsoft.com/en-us/research/project/non-transitive-dice.
30. Okuno, E., & Fratin, L. (2014). *Biomechanics of the Human Body. Undergraduate Lecture Notes in Physics*. New York: Springer.
31. Olteţeanu, A. M. (2016). From simple machines to Eureka in four not-so-easy steps: Towards creative visuospatial intelligence. In V. Müller (Ed), *Fundamental Issues of Artificial Intelligence* (pp 161-182). Synthese Library (Studies in Epistemology, Logic, Methodology, and Philosophy of Science), vol 376. Cham: Springer.
32. Pegg, E., Jr. (2005). Tournament dice. Retrieved from http://www.mathpuzzle.com/MAA/39-Tournament%20Dice/mathgames_07_11_05.html.
33. Poddiakov, A. (2010). Intransitivity cycles, and complex problem solving. Paper presented at the 2nd mini-conference "Rationality, Behavior, Experiment"; Moscow, September 1-3, 2010. https://www.researchgate.net/publication/237088961
34. Poddiakov, A. (2017). [Rule of transitivity *vs* intransitivity of choice]. *Nauka i Zhizn', 3*, 130-137. https://www.nkj.ru/archive/articles/30869/ (in Russian)
35. Poddiakov, A., & Valsiner, J. (2013). Intransitivity cycles and their transformations: How dynamically adapting systems function. In L. Rudolph (Ed.), *Qualitative Mathematics for The Social Sciences: Mathematical Models for Research on Cultural Dynamics* / Ed. by Abingdon, NY: Routledge. 343-391. Retrieved from https://www.researchgate.net/publication/281288415
36. Reichenbach, T., Mobilia, M., & Frey, E. (2007). Mobility promotes and jeopardizes biodiversity in rock–paper–scissors games. *Nature. 448,* 1046-1049. http://dx.doi.org/10.1038/nature06095.
37. Roberts, T. S. (2004). A ham sandwich is better than nothing: Some thoughts about transitivity. *Australian Senior Mathematics Journal, 18 (2),* 60–64. Retrieved from https://files.eric.ed.gov/fulltext/EJ720054.pdf.
38. Rogers, K. (2005). *On the Metaphysics of Experimental Physics*. London: Palgrave Macmillan.
39. Scheinerman, E.A. (2012). *Mathematics: A Discrete Introduction*. Belmont, CA: Brooks Cole.
40. Sinervo, B., & Lively, C. M. (1996). The rock-paper-scissors game and the evolution of alternative male strategies. *Nature*, *380*, 240–243. http://dx.doi.org/10.1038/380240a0
41. Strogatz, S. H. (2015). *Nonlinear dynamics and chaos: with applications to physics, biology, chemistry, and engineering*. Boulder, CO: Westview Press, a member of the Perseus Books Group.
42. Szolnoki, A., Mobilia, M., Jiang, L. L., Szczesny, B., Rucklidge, A.M., & Perc, M. (2014). Cyclic dominance in evolutionary games: A review. *Journal of the Royal Society, Interface*, *11*(100), 1-20. http://dx.doi.org/10.1098/rsif.2014.0735.

43. Tang, C.-Y, O'Brien, M. J., & Hawkins, G. F. (2005). Embedding simple machines to add novel dynamic functions to composites. *The Journal of The Minerals, Metals & Materials Society*, *57*(3), 2–35. https://doi.org/10.1007/s11837-005-0230-y.
44. Temkin, L. (1996). A continuum argument for intransitivity. *Philosophy & Public Affairs*, *25*, 175–210. https://doi.org/10.1111/j.1088-4963.1996.tb00039.x.
45. Temkin, L. (2012). *Rethinking the Good: Moral Ideals and the Nature of Practical Reasoning*. Oxford University Press.
46. Transitive law. Retrieved from https://www.britannica.com/topic/transitive-law.
47. Trybuła, S. (1961). On the paradox of three random variables. *Applicationes Mathematicae, 5*(4), 321-332. Retrieved from http://eudml.org/doc/264121.
48. Vandermeer, J. (2011). Intransitive loops in ecosystem models: From stable foci to heteroclinic cycles. *Ecological Complexity*, *8*(1), 92-97. https://doi.org/10.1016/j.ecocom.2010.08.001.
49. Wakild, T. (2006). Simple Machines. Grade 5. Retrieved from http://wmich.edu/engineer/ceee/edcsl/pdf/Simple%20Machine%20STEM%20guide.pdf.
50. West, L.J., & Hankin, R. (2008). Exact tests for two-way contingency tables with structural zeros. *Journal of Statistical Software, 28*(11), 1-19. http://dx.doi.org/10.18637/jss.v028.i11.
51. Williams, R. L., Chen, M. Y., & Seaton, J. M. (2003). Haptics-augmented simple-machine educational tools. *Journal of Science Education and Technology*, *12*(1), 1-12. https://doi.org/10.1023/A:1022114409119.